\documentclass{amsart}

\usepackage{amsmath, graphicx}
\usepackage{amsfonts}
\usepackage{amssymb}

\numberwithin{equation}{section}

\newtheorem{theorem}[subsection]{Theorem}
\newtheorem{proposition}[subsection]{Proposition}
\newtheorem{lemma}[subsection]{Lemma}

\newtheorem{conjecture}[subsection]{Conjecture}

\newcommand{\sumstar}{\mathop{{\sum}^{*}}}
\newcommand{\half}{\frac{1}{2}}
\newcommand{\sym}{\mathrm{sym}^2}

\newcommand{\intt}{\int_{-\infty}^{\infty}}

\newcommand{\ep}{\epsilon}
\newcommand{\hhalf}{\tfrac{1}{2}}

\theoremstyle{plain}

\begin{document}

\title{The average of $L^4$-norms of holomorphic cusp forms}
\author{Rizwanur Khan}
\address{Mathematisches Institut, Georg-August Universit\"{a}t G\"{o}ttingen, Bunsenstrasse 3-5, D-37073 G\"{o}ttingen, Germany}
\email{rrkhan@uni-math.gwdg.de}
\thanks{The author was supported by a grant from the European Research Council (grant agreement number 258713)}

 \begin{abstract}
We find the average value of the $L^4$-norm of holomorphic Hecke cusp forms of weight within a dyadic interval, up to an error which saves a power of the weight.
 \end{abstract}

\maketitle

\section{Introduction}

Let $f$ be a holomorphic cusp form of weight $k$ for $SL_2(\mathbb{Z})$ which is an eigenform of every Hecke operator. Define the rescaled function $F(z)= y^{k/2}f(z)$, for $z=x+iy$ in the upper half plane $\mathbb{H}$. The absolute value of $F$ is  $SL_2(\mathbb{Z})$-invariant, so we may ask how $|F|$ is distributed in terms of $k$ on the modular curve $SL_2(\mathbb{Z}) \backslash \mathbb{H}$. This interesting question has received a lot of attention recently with Holowinsky and Soundararajan's proof \cite{holsou} of a conjecture of Rudnick and Sarnak closely related to the Quantum Unique Ergodicity conjecture. It was proven that the $L^2$-mass of $F$ is equidistributed over $SL_2(\mathbb{Z}) \backslash \mathbb{H}$ for large $k$. Indeed, understanding the $L^p$-norms of $F$ is key to learning about its size distribution. For $1\le  p < \infty$, the $L^p$-norm of $F$ is defined to be 
\begin{align}
\| F \|_p = \left( \int_{SL_2(\mathbb{Z}) \backslash \mathbb{H}} |F(z)|^p \frac{dx dy}{y^2} \right)^{1/p}.
\end{align}

The $L^4$-norm was studied by Blomer, Khan, and Young in \cite{blokhayou}. They obtained the first non-trivial upper bound for $\| F \|_4$, proving with the normalization $\| F \|_2=1$ that
\begin{align}
\label{bkybound} \| F \|_4^4 \ll_{\ep}  k^{1/3+\epsilon}
\end{align}
for any $\epsilon>0$. However this seems to be very far from the truth. The following conjecture was made in the same paper. Note that 
\begin{align}
\int_{SL_2(\mathbb{Z}) \backslash \mathbb{H}}  \frac{dx dy}{y^2}=\frac{\pi}{3}.
\end{align}
\begin{conjecture}\cite[Conjecture 1.2]{blokhayou} \label{conj}
With the normalization  
\begin{align}
\label{normaliz} \frac{3}{\pi} \| F \|_2^2 =  1,
 \end{align}
we have that
\begin{align}
\label{conjec}  \frac{3}{\pi} \| F \|_4^4 \sim 2
\end{align}
as $k\to \infty$.
\end{conjecture}
\noindent Although this conjecture appears to be out of reach of current methods, it is a valuable step towards understanding what to expect for the distribution of $|F|$. In the parallel program of speculating on the distribution of a Hecke Maass cusp form $\phi$ of large eigenvalue for $SL_2(\mathbb{Z})$, predictions are already well established. The Random Wave conjecture, conceived by Berry \cite{ber} and investigated further by Hejhal and Rackner \cite{hejrac}, says that $\phi$ should obey the normal distribution law. Thus the conjecture analogous to (\ref{conjec}) is that
\begin{align}
\frac{3}{\pi} \| \phi \|_4^4 \sim 3
\end{align}
with the normalization $\frac{3}{\pi} \| \phi \|_2^2 =  1$, which coincides with the the fourth moment of a standard normal random variable.  

Our goal is to establish Conjecture \ref{conj} on average over $f$ and $k$. Note that we use the normalization $\| F \|_2 =1$ instead of (\ref{normaliz}).
\begin{theorem}\label{thm}  Let $w$ be a smooth, non-negative function supported on $(1,2)$ with bounded derivatives. Let $B_k$ be the orthonormal Hecke basis of the space of holomorphic cusp forms of weight $k$. Then for some $\delta>0$ we have that
\begin{align}
 \frac{2}{KW} \sum_{k \equiv 0 \bmod 2} w\Big(\frac{k}{K}\Big) \frac{12}{k} \sum_{f \in B_k} \left\| F  \right\|_4^4 = \frac{6}{\pi} + O(K^{-\delta}),
\end{align}
where $W=\int_{-\infty}^{\infty} w(x) dx.$
In particular, $\| F \|_4^4\ll 1$ for almost all $f$ of weight $K\le k \le 2K$. 
\end{theorem}

The main steps of the proof are as follows. By a formula of Watson, the theorem is reduced to finding the mean value of a degree 8 $L$-function at the central point. We must find the triple product $L$-function value\begin{align}
\label{tripprod} L(\hhalf, f\times f\times g) = L(\hhalf,\sym f \times g)L(\hhalf, g)
\end{align}
on average over $g\in B_{2k}$, $f\in B_k$, and $k$ of size $K$, which amounts to the average of about $K^3$ terms over about $K^3$ harmonics, by its approximate functional equation. The Petersson trace formula leads to the task of bounding a sum like
\begin{align}
\label{twoj} \sum_{n,m,c_1,c_2} S(n,m;c_1) S(n^2, 1;c_2)  \sum_{\substack{k\sim K\\ 2|k}}  i^k J_{2k-1}\Big( \frac{4\pi \sqrt{nm}}{c_1} \Big) J_{k-1}\Big(\frac{4\pi n}{c_2} \Big) ,
\end{align}
where $n,m,c_1,c_2$ are positive integers whose range of summation depends on $K$. The interesting inner sum was studied by Liu and Young in \cite{liuyou}, but unfortunately their analysis falls short of our needs. In that paper, the $L^2$-norm of certain Siegel modular forms is investigated through the mean value of a degree 6 $L$-function similar to the one appearing on the right hand side of (\ref{tripprod}). The extra $GL(2)$ factor in this line makes our problem more complex. Building on the beautiful work of Liu and Young, we are able to improve their treatment of the inner sum in (\ref{twoj}). It turns out that the phase of this sum is a rational number, which allows us to perform Poisson summation in $n$ and carefully estimate certain exponential sums over finite fields at the end of the proof.

{\bf Notation.} Throughout, $\epsilon$ will denote an arbitrarily small positive constant, but not necessarily the same one from one occurrence to the next. In any estimate, any implied constant may depend on $\epsilon$. The delta symbol $\delta_{P}$ will equal 1 whenever the statement $P$ is true and 0 whenever it is false.

\subsection{$L$-functions}

By defintion, $f\in B_k$ is an eigenform of every Hecke operator $T_n$. Let
\begin{align}
\label{eigenvalues} T_n f = a_f(n)n^{\frac{k-1}{2}} f.
\end{align}
The rescaled eigenvalues $a_f(n)$ satisfy the Hecke relations
\begin{align}
\label{heckerel} a_f(n)a_f(m) = \sum_{d|(n,m)}a_f\left(\frac{nm}{d^2}\right),
\end{align}
Deligne's bound
\begin{align}
a_f(n)\le \tau(n),
\end{align}
where $\tau(n)$ is the number of divisors of $n$, and the following orthogonality relation.
\begin{lemma}{\bf (Petersson trace formula)}
We have
\begin{align}
\frac{2\pi^2}{k-1} \sum_{f\in B_k} \frac{a_f(n)a_f(m)}{L(1, \sym f)} = \delta_{n=m} +2\pi i^{-k} \sum_{c\ge 1} \frac{S(n,m;c)}{c}J_{k-1}\left(\frac{4\pi \sqrt{nm}}{c}\right),
\end{align}
where $L(s,\sym f)$ is the symmetric-square $L$-function defined in (\ref{symsq}), $J_{k-1}$ is the Bessel function defined in (\ref{jdef}) and 
\begin{align}
S(n,m;c) = \sumstar_{b \bmod c} e\left( \frac{nb + m\overline{b}}{c} \right)
\end{align}
is the Kloosterman sum.
\end{lemma}
\proof
See \cite[Proposition 14.5]{iwakow}.
\endproof
\noindent The Kloosterman sum satisfies the multiplicative property
\begin{align}
\label{kloosmult} S(n,m,c_1c_2)= S(n,m\overline{c_2}^2;c_1)S(n,m\overline{c_1}^2;c_2)
\end{align}
for $c_1$ and $c_2$ coprime and Weil's bound 
\begin{align}
\label{klobound} |S(n,m;c)| \le \tau(c)c^{\half}(n,m,c)^{\half}.
\end{align}
Note that on average over $n,m$ or $c$, this bound is $O((nmc)^{\ep}c^{1/2})$. For example we have that
\begin{align}
\sum_{n<x}  |S(n,m;c)| \ll c^{\half+\ep }\sum_{n<x} (n,c)^{1/2} \ll c^{\half+\ep } \sum_{d|c} d^{1/2} \sum_{\substack{n<x\\ d|n}} 1  \ll c^{\half+\ep }x.
\end{align}

We will work with the following entire $L$-functions, defined for for $f\in B_k, g\in B_{2k}$ and $\Re(s) > 1$ by
\begin{align}
&L(s,f) = \sum_{n\ge 1} \frac{a_f(n)}{n^s},\\
\label{symsq} &L(s, \sym f)= \sum_{n\ge 1} \frac{A_f(n,1)}{n^s} = \zeta(2s) \sum_{n\ge 1} \frac{a_f(n^2)}{n^s},\\
&L(s, \sym f \times g) = \sum_{n,r\ge 1} \frac{A_f(n,r)a_g(n)}{(nr^2)^s} = \sum_{n,r,d\ge 1}\mu(d) \frac{A_f(n,1)A_f(r,1)a_g(nd)}{(nr^2d^3)^s},
\end{align}
where $A_f(n,r)=A_f(r,n)$, defined implicitly, are the Fourier coefficients of the symmetric square lift of $f$, a cusp form on $GL(3)$. By the normalization (\ref{eigenvalues}), these $L$-functions have functional equations relating values at $s$ and $1-s$, so that the central point is always $s=1/2$.

We observe as in \cite{blokhayou} the decomposition
\begin{align}
\| F \|_4^4 = \langle F^2 , F^2 \rangle = \sum_{g\in B_{2k}} | \langle F^2 , G  \rangle |^2,
\end{align}
which holds because $f^2$ is a cusp form of weight $2k$.
 By Watson's formula \cite[Theorem 3]{wat} we have 
 \begin{align}
 |\langle F^2, G\rangle |^2 =  \frac{\pi^3}{2(2k-1)}  \frac{ L(s, g)  L(s, \text{sym}^2f  \times g)  }{L(1, \text{sym}^2 f)^2 L(1, \text{sym}^2 g)} . 
\end{align}
Thus
\begin{align}
\label{watson} \| F \|_4^4 = \frac{\pi^3}{2(2k-1)L(1, \text{sym}^2 f)^2}  \sum_{g \in B_{2k}} \frac{ L(\half, g)  L(\half, \text{sym}^2f  \times g)  }{  L(1, \text{sym}^2 g)}.
\end{align}

We get a handle on the central values of $L$-functions using approximate functional equations.
\begin{lemma}{\bf (Approximate functional equation)}
Let
\begin{align}
\Lambda_{k,1}(s)= (2\pi)^{-s} \Gamma(s+k-\tfrac{1}{2})
\end{align}
and
\begin{align}
\Lambda_{k,2}(s)= 8(2\pi)^{-3s-3k+\frac{3}{2}}\textstyle\Gamma(s + 2k - \frac{3}{2})\Gamma(s+ k - \frac{1}{2})\Gamma(s+\frac{1}{2}).
\end{align}
For $\xi,\sigma >0$ and $j=1,2$, define
\begin{align}
V_{k,j}(\xi)=\frac{1}{2\pi i} \int_{(\sigma)}\frac{\Lambda_{k,j}(1/2+s)}{\Lambda_{k,j}(1/2)} \xi^{-s} \frac{ds}{s}.
\end{align}
For $g \in B_{2k}$ we have
\begin{align}
\label{sum1} L(1/2,g)= 2\sum_{m\ge 1} \frac{a_g(m)}{m^{1/2}}V_{k,1}(m)
\end{align}
and for $f\in B_{k}$ we have
\begin{align}
\label{sum2} L(1/2,\sym f\times g)= 2\sum_{n,r\ge 1} \frac{A_f(n,r)a_g(n)}{(nr^2)^{1/2}}V_{k,2}( n r^2).
\end{align}
\end{lemma}
\proof
This follows from \cite[Theorem 5.3]{iwakow} and the functional equations
\begin{align}
\Lambda_{k,1}(s)L(s,g)=\Lambda_{k,1}(1-s)L(1-s,g)
\end{align}
and
\begin{align}
\Lambda_{k,2}(s)L(s, \sym f \times g)=\Lambda_{k,2}(1-s)L(1-s, \sym f \times g),
\end{align}
which may be found in \cite[Chapter 7]{iwa} and \cite[Section 4.1]{wat} respectively.
\endproof
\noindent By Stirling's formula, we have for $\Re(s)=A$ that
\begin{align}
\label{lambdabound} \frac{\Lambda_{k,j}(1/2+s)}{\Lambda_{k,j}(1/2)} \ll_A k^{j A},
\end{align}
except for where the left hand side has a pole.
It follows that for any $A>0$ and integer $B\ge 0$ we have
\begin{align}
\label{vbound} \Big(\frac{\xi}{k^j} \Big)^{B} V_{k,j}^{(B)}(\xi) \ll_{A,B} \Big(1+\frac{\xi}{k^j} \Big)^{-A}.
\end{align}
Thus the sums in (\ref{sum1}) and (\ref{sum2}) are essentially supported on $m<k^{1+\epsilon}$ and $nr^2<k^{2+\epsilon}$ respectively. So in Theorem \ref{thm} we are averaging about $K^{3}$ terms over about $K^3$ harmonics, as $|B_k|\sim k/12$. We also note that for $u\in (1,2)$ we have
\begin{align}
\label{v-h} \frac{d^B}{du^B} V_{uK,j} (\xi) \ll_{A,B} \Big(1+\frac{\xi}{K^j} \Big)^{-A},
\end{align}
by Stirling's formula.

For the $L$-function values at the edge of the critical strip, we will need
\begin{lemma} \label{remove}
Given $0<\delta<1/10$, we have for all but $O(k^{100\delta})$ forms in $B_k$ that
\begin{align}
\label{1series} L(1,\sym f)^{-1} =  \sum_{d_1,d_2,d_3\ge 1} \frac{\mu(d_1d_2d_3)\mu(d_2) a_f(d_1^2d_2^2)}{d_1d_2^2d_3^3}\exp\left(-\frac{d_1d_2^2d_3^3}{k^{\delta}}\right) + O(k^{-\delta^2+\epsilon}).
\end{align}
\end{lemma}

\proof
Let $0<\delta<1/10$. By a zero-density result of Lau and Wu \cite[(1.11)]{lauwu}, for all but $O(k^{100\delta})$ forms $f$ in $B_k$ we have that $L(s,\sym f)$ is nonzero in the region
\begin{align}
\mathcal{R}= \{ s: \Re(s) \ge 1- 3\delta, |\Im(s)|\le  3k^{\delta} \}.
\end{align}
Thus the lemma is proved if we can establish (\ref{1series}) on the assumption that $L(s,\sym f)$ has no zero in $\mathcal{R}$.

We first show that
\begin{align}
\label{l'lclaim0} L(s,\sym f)^{-1} \ll k^{\epsilon}
\end{align}
for $s=\sigma +it \in \mathcal{R}'=  \{ s: \Re(s) \ge 1- \delta, |\Im(s)|\le  k^{\delta} \}$. We have that
\begin{align}
\log L(s,\sym f) = \log L(2+it) - \int_{\sigma}^{2} \frac{L'}{L} (x+it) dx.
\end{align}
Thus the claim would follow from showing that 
\begin{align}
\label{l'lclaim} \frac{L'}{L}(s, \sym f) =  o(\log k )
\end{align}
for $s \in \mathcal{R}'$. We follow the ideas in \cite[Proposition 5.17]{iwakow}. By Cauchy's theorem, the exponential decay of the gamma function in vertical lines and the bound $ \frac{L'}{L}(s, \sym f) \ll \log k$ of \cite[Theorem 5.7]{iwakow} which holds for all $f\in B_k$, we have for $1<X< k$ that
\begin{multline}
\label{l'l} \frac{L'}{L}(s, \sym f) =\\ \frac{1}{2\pi i}  \int_{2-ik^{\delta}}^{2+ik^{\delta}} \frac{L'}{L}(w+s,\sym f) \Gamma(w) X^{w} dw -  \frac{1}{2\pi i}\int_{-\delta-ik^{\delta}}^{-\delta+ik^{\delta}} \frac{L'}{L}(w+s,\sym f) \Gamma(w) X^{w} dw + O(k^{-100}).
\end{multline}
By the Euler product
\begin{align}
L(w, \sym f) = \prod_p (1-a_f(p^2)p^{-w} + a_f(p^2)p^{-2w} - p^{-3w})^{-1},
\end{align}
valid for $\Re(w)>1$, we have that 
\begin{align}
\frac{1}{2\pi i}  \int_{2-ik^{\delta}}^{2+ik^{\delta}} \frac{L'}{L}(w+s,\sym f) \Gamma(w) X^{w} dw = \sum_{n\ge 1} \frac{\Lambda_f(n)}{n^s} e^{-n/X} + O(k^{-100}),
\end{align}
where $\Lambda_f(n)\ll \log n$. Thus this integral is bounded by $X^{\delta+\epsilon}$. The integral on the line $\Re(w)=-\delta$ in (\ref{l'l}) is bounded by $X^{-\delta}\log k$. Taking $X=(\log k)^{1/2}$, say, proves the claim (\ref{l'lclaim}).

Now by Cauchy's theorem again we have
\begin{align}
L(1,\sym f)^{-1} = \frac{1}{2\pi i} \int_{2-ik^{\delta}}^{2+ik^{\delta}} \frac{ \Gamma(s) k^{\delta s}}{L(1+s,\sym f)} ds -  \frac{1}{2\pi i}\int_{-\delta-ik^{\delta}}^{-\delta+ik^{\delta}} \frac{\Gamma(s) k^{\delta s}}{L(s, \sym f)}  ds + O(k^{-100}).
\end{align}
The first integral equals
\begin{align}
\frac{1}{2\pi i}\int_{2-ik^{\delta}}^{2+ik^{\delta}} \frac{ \Gamma(s) k^{\delta s}}{L(1+s,\sym f)} ds = \sum_{d_1,d_2,d_3\ge 1} \frac{\mu(d_1d_2d_3)\mu(d_2) a_f(d_1^2d_2^2)}{d_1d_2^2d_3^3}\exp\left(-\frac{d_1d_2^2d_3^3}{k^{\delta}}\right) + O(k^{-100}).
\end{align}
For the integral inside the critical strip, we have by (\ref{l'lclaim0}) that
\begin{align}
  \int_{-\delta-ik^{\delta}}^{-\delta+ik^{\delta}} \frac{\Gamma(s) k^{\delta s}}{L(1+s, \sym f)}  ds\ll k^{-\delta^2+\epsilon}.
 \end{align}
 \endproof

By Lemma \ref{remove}, the bounds
\begin{align}
k^{-\epsilon}\ll L(1,\sym f) \ll k^{\epsilon}
\end{align}
and the trivial bound $\left\| F  \right\|_4^4\ll k^{1/4+\ep}$, we see that Theorem \ref{thm} would follow from showing that there exists some $\delta>0$ such that for any $1\le d_1,d_2,d_3\le k^{\ep}$, we have
\begin{align}
\label{ffrom} \frac{2}{K W} \sum_{k \equiv 0 \bmod 2} w\left(\frac{k}{K}\right) \frac{12}{k} \sum_{f \in B_k}L(1,\sym f) a_f(d_1^2d_2^2) \left\| F  \right\|_4^4 = \frac{6}{\pi} + O(K^{-\delta }).
\end{align}

\subsection{Bessel functions}
For $\ell>0$ an integer, the $J$-Bessel function
\begin{align}
\label{jdef} J_\ell (x) = \int_{-\half}^{\half} e(\ell t)  e^{-i x \sin (2\pi t)} \ dt
\end{align}
satisfies the bounds (cf. \cite[Lemma 4.2-4.3]{ran})
\begin{align} \label{jbound}
J_{\ell}(x) \ll
\begin{cases}
e^{-\ell} &\text{ for } x< \ell/10,\\
\min(\ell^{-1/3}, |x^2 - \ell^2|^{-1/4}) &\text{ for } x\ge \ell/10.
\end{cases}
\end{align}
It may be approximated using Langer's formula below. 
\begin{lemma} \label{langer}
For $x> \ell>0$, let
\begin{align}
w =  \left( \frac{x^2}{\ell^2} - 1\right)^{\frac{1}{2}} 
\end{align}
and
\begin{align}
z = \ell\left(w-\tan^{-1} w \right).
\end{align}
We have
\begin{align}
J_\ell (x) = \left(1-\frac{\tan^{-1} w}{w}\right)^{1/2} \left(\frac{2}{\pi z} \right)^{1/2}  \bigg( \cos\left( z -\frac{\pi}{4} \right) + O\left( \frac{1}{z} \right)\bigg) + O\left(\frac{1}{\ell^{4/3}}\right).
\end{align}
\end{lemma}
\proof
See \cite[Sections 7.13.1,7.13.4]{erd}.
\endproof
The value of $J_{k-1}(y)$ on average over integers $k$ divisible by 4 is well understood. We record only an upper bound.
\begin{lemma}  \label{noik}
Let $0<y<K^{2+\ep}$ and $h$ a smooth, non-negative function supported on $(1,2)$ such that
\begin{align}
\label{h} h^{(B)}(u) \ll_B K^{\ep B}
\end{align}
for any integer $B\ge 0$. We have that
\begin{align}
\label{agn} 4\sum_{k \equiv 0 \bmod 4} h\Big(\frac{k}{K}\Big) J_{k-1}( y) \ll h\Big(\frac{y}{K}\Big)+\frac{1}{K^{1-\epsilon}} + \frac{y}{K^{3-\epsilon}}.
\end{align}
\end{lemma}
\proof
See \cite[Lemma 5.8]{iwa}. In the result of that lemma, the first main term equals $h(y/K)$, the second main term is bounded by $K^{-1+\ep}$, and the error term is less than $yK^{-3+\ep}$.
\endproof

We also need to understand the following average value of a product of two Bessel functions:
\begin{align}
\label{prodbes} \sum_{k \equiv 0 \bmod 2} i^k h\Big(\frac{k}{K}\Big) J_{k-1}( x) J_{2k-1}(y)
\end{align}
for $0<x<K^{2+\ep}$ and $0<y<K^{3/2+\ep}$. A very similar sum was found by Liu and Young in \cite[Theorem 5.3]{liuyou} for $x<K^{2+\ep}$ and $y<K^{1+\ep}$, up to an error of $O(K^{-1+\ep})$. We must extend the range of admissible $y$ and improve the error term of their asymptotic. The next result achieves this for $x<K^{2-\ep}$.

\begin{lemma} \label{doubleavg}
Let $A>0$ and $h$ a smooth, non-negative function supported on $(1,2)$ with bounded derivatives. Let $x$ and $y$ belong to the set
\begin{align}
\{0< x \le K^{4/3-\ep}, y>0  \} \cup \{ K^{4/3-\ep} < x < K^{2-\ep} , 0< y< xK^{-\ep}  \} \cup  \{ K^{4/3-\ep} < x <K^{2-\ep}, y> xK^{\ep}  \}.
\end{align}
There exists smooth functions
$H^\pm(\xi_1;\xi_2,\xi_3,\xi_4)$ depending on $A$ with bounded derivatives in any compact set, polynomial in $\xi_2,\xi_3,\xi_4$ and supported on
\begin{align}
\label{pro1} 1 \ll |\xi_1| \ll 1,
\end{align}
such that
\begin{multline}
\label{avg2}  \sum_{k \equiv 0 \bmod 2} i^k h\Big(\frac{k}{K}\Big) J_{k-1}( 4\pi x) J_{2k-1}(4\pi y) = \\
\frac{1}{\sqrt{x}} \sum_{\pm} e\left( \pm \left(  \frac{y^2}{4x} +  2x\right) \right)  \Big(1\pm \frac{iy}{4x\sqrt{1-(\frac{y}{4x})^2}}\Big) H^{\pm}\Big(\frac{ y\sqrt{1- (\frac{y}{4x})^2}}{2K}\Big) + O_A(K^{-A}),
\end{multline}
where 
\begin{align}
\label{Hpm} H^{\pm}\Big(\frac{ y\sqrt{1- (\frac{y}{4x})^2}}{2K}\Big) = H^{\pm}\Big(\frac{ y\sqrt{1- (\frac{y}{4x})^2}}{2K}; \frac{x(1-2(\frac{y}{4x})^2)}{K^2},\frac{ y\sqrt{1-(\frac{y}{4x})^2}}{2K^3}, \frac{ y}{4x}\Big)\delta_{\frac{y}{4x}<1}.
\end{align}

Property (\ref{pro1}) implies that $H^{\pm}\Big(\frac{ y\sqrt{1- (\frac{y}{4x})^2}}{2K}\Big)$ is nonzero only for $x,y\gg K$ and supported on
\begin{align}
\label{pro2} K^2/y^2\ll 1-y/4x \ll K^2/y^2.
\end{align}
\end{lemma}

\noindent {\bf Remark.} The asymptotic (\ref{avg2}) cannot hold for $x>K^{4/3+\ep}$ and $y\sim x$ satisfying (\ref{pro2}). This is because the main term on the right hand side of (\ref{avg2}) would have size about $\sqrt{x}/K$ while the left hand side of (\ref{avg2}) is less than a constant multiple of  $K/x$, by (\ref{jbound}).

\proof
The proof is the same as that of \cite[Theorem 5.3]{liuyou} until the appearance of the oscillatory integral (\ref{uint2}), which we evaluate somewhat differently.

By (\ref{jdef}) we have that
\begin{multline}
 \sum_{k \equiv 0 \bmod 2} i^k h\Big(\frac{k}{K}\Big) J_{k-1}(4\pi  x) J_{2k-1}(4\pi y) \\
= \int_{-\half}^{\half}  \int_{-\half}^{\half}  e(-2x\sin( 2\pi t)-2y\sin(2\pi u)-t-u)   \sum_{k \equiv 0 \bmod 2} i^k h\Big(\frac{k}{K}\Big) e(k(t+2u))  \ du dt .
\end{multline}
By Poisson summation, this equals
\begin{align}
&  \frac{K}{2}\int_{-\half}^{\half}  \int_{-\half}^{\half}  e(-2x\sin( 2\pi t)-2y\sin(2\pi u)-t-u) \sum_{l \in \mathbb{Z}} \hat{h}\left(K\left(t+2u+\frac{l}{2}+\frac{1}{4}\right)\right) \ du dt \\
\label{tint} =&   \sum_{\pm} \frac{\mp i}{2} \int_{-\infty}^{\infty}  \hat{h}(t) \int_{-\half}^{\half} e\left( \pm 2x\cos\left( \frac{2\pi t}{K}-4\pi u\right)-2y\sin(2\pi u)-\frac{t}{K}+u\right) \ du dt,
\end{align}
where
\begin{align}
\hat{h}(t) = \intt h(v)e(-vt)\ dv
\end{align}
denotes the Fourier transform of $h$.
Since
\begin{align}
 \hat{h}(t) \ll_B (1+|t|)^{-B}
\end{align}
for any $B>0$ by repeated integration by parts, we may restrict the $t$-integral in (\ref{tint}) to $|t|\ll K^{\ep}$. By Taylor's theorem, we have for any $A>0$ that
\begin{align}
&e\left( \pm 2x\cos\left( \frac{2\pi t}{K}-4\pi u\right)\right)\\
\nonumber =&e\left( \pm 2x \cos  \frac{2\pi t}{K}\cos 4\pi u \pm 2x   \sin \frac{2\pi t}{K}\sin 4\pi u \right)\\
\nonumber = &e\left(  \pm 2x \cos 4\pi u \pm 4\pi tx K^{-1}\sin 4\pi u \mp 4\pi^2 t^2 xK^{-2} \cos 4\pi u \mp \tfrac{8}{3} \pi^3 t^3 x K^{-3} \sin 4\pi u + \ldots O_A( K^{-A})\right),
\end{align}
where the ellipsis represents finitely many terms of the Taylor expansion.
Further, since $x<K^{2-\epsilon}$, we have that
\begin{align}
e(\mp 4\pi^2 t^2 xK^{-2} \cos 4\pi u)= 1 \mp 4\pi^2 t^2 xK^{-2} \cos 4\pi u  + 8(\pi^2 t^2 xK^{-2} \cos 4\pi u)^2 + \ldots + O_A(K^{-A}),
\end{align}
\begin{align}
e(\mp \tfrac{8}{3} \pi^3 t^3 x K^{-3} \sin 4\pi u ) = 1 \mp \tfrac{8}{3} \pi^3 t^3 x K^{-3} \sin 4\pi u + \tfrac{32}{9}( \pi^3 t^3 x K^{-3} \sin 4\pi u)^2 +\ldots + O_A( K^{-A}),
\end{align}
and so on.
Thus we have that (\ref{tint}) equals
\begin{multline}
 \label{tint2} \sum_{\pm} \frac{\mp i}{2}  \int_{-\half}^{\half} e(u \pm 2x \cos 4\pi u-2y\sin 2\pi u )\int_{-\infty}^{\infty}  \hat{h}(t) e\left(t \left(\frac{-1 \pm 4\pi x \sin 4\pi u }{K} \right)\right) \\
 \times \left(1\mp \frac{8 \pi^3 t^3 x \sin 4\pi u }{3 K^3}+\ldots \right) \left(1 \mp \frac{4\pi^2 t^2 x \cos 4\pi u}{K^2}+\ldots\right)\cdots \ dt du + O_A(K^{-A}),
\end{multline}
Now by the identity
\begin{align}
\label{prope} \int_{-\infty}^{\infty}   \hat{h}(t) (2\pi it)^j e(v t) \ dt =h^{(j)}(v),
\end{align}
we have that the main term of (\ref{tint2}) equals
\begin{align}
\label{plusminus} &\sum_{\pm}   \int_{-\half}^{\half} e(u \pm 2x \cos 4\pi u-2y\sin 2\pi u  )  H_1^\pm \Big( \frac{ x \sin 4\pi u }{K}; \frac{ x \cos 4\pi u }{K^2}, \frac{ x \sin 4\pi u }{K^3}\Big) \ du \\
\label{uint}  =&\sum_{\pm}   \int_{-\frac{1}{4}}^{\frac{1}{4}} e(u \pm 2x \cos 4\pi u-2y\sin 2\pi u  )  H_2^\pm \Big( \frac{ x \sin 4\pi u }{K}; \frac{ x \cos 4\pi u }{K^2}, \frac{ x \sin 4\pi u }{K^3}\Big) \ du
\end{align}
for some functions $H_1^\pm(\xi_1; \xi_2,\xi_3)$ and $H_2^\pm(\xi_1; \xi_2,\xi_3)$ with bounded derivatives in any compact set, polynomial in $\xi_2,\xi_3$ and supported on
\begin{align}
\label{support} 1 \ll |\xi_1| \ll 1.
\end{align}

Making the substitution $v=\sin 2\pi u$, we get that (\ref{uint}) equals
\begin{align}
\label{uint2} \sum_{\pm}   \int_{-\infty}^{\infty} e (\pm 2x(1-2v^2)- 2yv ) \psi^\pm (v) \ dv,
\end{align}
where we define
\begin{align}
\psi^{\pm}(v) = 
\begin{cases}
\frac{i}{4\pi} \Big(1+\frac{iv}{\sqrt{1-v^2}}\Big) H_2^{\pm}\Big(\frac{2xv\sqrt{1-v^2}}{K}; \frac{x(1-2v^2)}{K^2},\frac{2xv\sqrt{1-v^2}}{K^3}\Big) \hfill \text{ for } |v|\le 1,\\
0 \hfill \text{ for} |v|> 1.
\end{cases}
\end{align}
By (\ref{support}), we have that $\psi^{\pm}(v)$ vanishes unless $0<|v|<1$ and $x\gg K$, in which case it is supported on 
\begin{align}
\label{supp0} \frac{K}{x}\ll |v| \ll \frac{K}{x}
\end{align}
and
\begin{align}
\label{supp1} 1-\frac{K^2}{x^2} \ll |v| \ll 1-\frac{K^2}{x^2}.
\end{align}
For bounds, we have that
\begin{align}
\label{Hhat} \frac{\partial^B}{\partial v^B}  \psi^{\pm}(v)  \ll_B
 \begin{cases}
 \big(\frac{x}{K}\big)^{B} &\text{ if } |v|<1/2,\\
  \big(\frac{x}{K}\big)^{2B+1} &\text{ if } |v|\ge 1/2
 \end{cases}
\end{align}
 for any integer $B\ge 0$.
  
From this point, our proof differs from \cite{liuyou}. On observing that
\begin{align}
 e(\pm 2x(1-2v^2)-2yv) = e\Big(\mp4x\Big( v \pm \frac{y}{4x} \Big)^2 \pm \frac{y^2}{4x} \pm 2x \Big),
\end{align}
we get that (\ref{uint2}) equals
\begin{align}
\label{until}  \sum_{\pm} e\Big(  \pm \frac{y^2}{4x} \pm 2x \Big) \int_{-\infty}^{\infty} \psi^{\pm}(v) e\Big(\mp4x\Big( v \pm \frac{y}{4x} \Big)^2 \Big) \ dv.
\end{align}
 {\bf Case 1:} Assume that $0<x\le K^{4/3-\ep}$ and $y>0$.
 By Parseval's theorem we have that
\begin{align}
 \label{parseval} \int_{-\infty}^{\infty} e\Big(\mp4x\Big( v \pm \frac{y}{4x} \Big)^2 \Big)\psi^{\pm}(v) \ dv =  \frac{e(\tfrac{1}{8})}{\sqrt{8x}} \int_{-\infty}^{\infty} e\Big( \pm \frac{ t^2}{16 x}\pm \frac{yt}{4x}\Big) \widehat{\psi^{\pm}}(t) \  dt.
 \end{align}
 By (\ref{Hhat}) we may restrict the integral to $|t|< \frac{x^2}{K^{2-\ep}}$, up to an error of $O_A(K^{-A})$. The condition $x\le K^{4/3-\epsilon}$ implies that
 \begin{align}
\label{t2x} \Big| \frac{t^2}{x} \Big| < K^{-\epsilon}.
\end{align}
Taking a Taylor series expansion of $e\big(\frac{t^2}{16x}\big)$ and using property (\ref{prope}), we get that (\ref{parseval}) equals, up to an error of $O_A(K^{-A})$,
\begin{align}
\label{add1}\frac{1}{\sqrt{x}} \Big(1\pm \frac{iy}{4x\sqrt{1-(\frac{y}{4x})^2}}\Big) H^{\pm}\Big(\frac{\pm y\sqrt{1- (\frac{y}{4x})^2}}{2K}; \frac{x(1-2(\frac{y}{4x})^2)}{K^2},\frac{\pm y\sqrt{1-(\frac{y}{4x})^2}}{2K^3}, \frac{\pm y}{4x}\Big)\delta_{y/4x<1},
\end{align}
for some functions $H^\pm(\xi_1; \xi_2,\xi_3, \xi_4)$ described in the statement of the lemma.

 {\bf Case 2:} Assume that $K^{4/3- \ep}<x<K^{2-\ep}$ and $y<xK^{-\ep}$ or $y>xK^{\ep}$. By (\ref{supp0}-\ref{supp1}) and the condition $x>K^{4/3-\ep}$, we have $\psi^\pm (v)$ is supported on disjoint sets of size $K/x$ near $0$ and size $K^2/x^2$ near $1$ and $-1$. Thus we may split the integrals in (\ref{until}). We have that
 \begin{align}
\label{split}  &\int_{-\infty}^{\infty} e\Big(\mp 4x\Big( v \pm \frac{y}{4x} \Big)^2 \Big) \psi^{\pm}(v) \ dv= \\
 \label{split2} &\int_{-\infty}^{\infty} e\Big(\mp4x\Big( v \pm \frac{y}{4x} \Big)^2 \Big) \psi^{\pm}(v) \delta_{|v|<1/2} \ dv  +  \int_{-\infty}^{\infty} e\Big(\mp4x\Big( v \pm \frac{y}{4x} \Big)^2 \Big) \psi^{\pm}(v) \delta_{|v|>1/2 } \ dv.
 \end{align}
To the first integral in (\ref{split2}) we apply Parseval's theorem just as in (\ref{parseval}). Since $x<K^{2-\ep}$ and the derivatives of $\psi^\pm (v)$ are not so large for $|v|\ll K/x$, by (\ref{Hhat}), we get that (\ref{t2x}) holds. Then just as in (\ref{add1}) we get that this integral equals, up to an error of $O_A(K^{-A})$,
 \begin{align}
\label{split1ans} \Big(1\pm \frac{iy}{4x\sqrt{1-(\frac{y}{4x})^2}}\Big) H^{\pm}\Big(\frac{\pm y\sqrt{1- (\frac{y}{4x})^2}}{2K}; \frac{x(1-2(\frac{y}{4x})^2)}{K^2},\frac{\pm y\sqrt{1-(\frac{y}{4x})^2}}{2K^3}, \frac{\pm y}{4x}\Big) \delta_{y/4x< 1/2}.
 \end{align}
If $y<xK^{-\ep}$ then $ \delta_{y/4x<1/2}=1$. If $y>xK^{\ep}$, then $\delta_{y/4x<1/2}=0$ and (\ref{Hpm}) is also zero. 

As for the second integral in (\ref{split2}), we integrate by parts $B$ times. By (\ref{Hhat}) and the assumptions on $y$, we get that this integral is bounded by
\begin{align}
 x^{-B} \Big(\frac{x}{K}\Big)^{2B+1}.
\end{align}
This is less than any negative power of $K$ for large enough $B$, since $x<K^{2-\ep}$.
\endproof

Lemma \ref{doubleavg} does not cover all the ranges of $x$ and $y$ that we need. In the remaining ranges, we find suitable upper bounds for (\ref{prodbes}).

\begin{lemma}\label{doublebound}
For $K^{4/3-\ep}<x<K^{2-\ep}$, $xK^{-\ep}< y <xK^{\ep}$ and $h$ as in Lemma \ref{doubleavg}, we have that
\begin{align}
\label{doublebound1} \sum_{k \equiv 0 \bmod 2} i^k h\Big(\frac{k}{K}\Big) J_{k-1}( 4\pi x) J_{2k-1}(4\pi y)\ll 
 \frac{K}{\sqrt{xy}}.
 \end{align}
 If further $x$ and $y$ satisfy
 \begin{align}
 \label{furth} \Big|1-\frac{y}{4x}\Big|> \frac{K^{2+\ep}}{x^2},
 \end{align} 
 then we have that
 \begin{align}
\label{doublebound2} \sum_{k \equiv 0 \bmod 2} i^k h\Big(\frac{k}{K}\Big) J_{k-1}( 4\pi x) J_{2k-1}(4\pi y)\ll_B
\frac{x}{K} \Big(\frac{x^2}{K^2|4x-y|}\Big)^{B} 
 \end{align}
 for any integer $B\ge 0$.
\end{lemma}
\proof
The first estimate (\ref{doublebound1}) follows directly from (\ref{jbound}). For the second estimate, our starting point is the inequality
\begin{align}
& \sum_{k \equiv 0 \bmod 2} i^k h\Big(\frac{k}{K}\Big) J_{k-1}( 4\pi x) J_{2k-1}(4\pi y) \le \\
 \label{split3} & \Big|\int_{-\infty}^{\infty} e\Big(\mp4x\Big( v \pm \frac{y}{4x} \Big)^2 \Big) \psi^{\pm}(v) \delta_{|v|<1/2} \ dv \Big| +  \Big|\int_{-\infty}^{\infty} e\Big(\mp4x\Big( v \pm \frac{y}{4x} \Big)^2 \Big) \psi^{\pm}(v) \delta_{|v|>1/2 } \ dv \Big|
\end{align}
which follows from (\ref{split2}). The first integral in (\ref{split3}) is given by (\ref{split1ans}), up to any negative power of $K$. But since $y>K^{4/3-\ep}$, we have by (\ref{pro2}) that (\ref{split1ans}) vanishes unless $y/4x\sim 1$. However in that case it would vanish anyway since then $\delta_{y/4x<1/2}=0$. The second integral in (\ref{split3}) can be integrated by parts $B$ times. By (\ref{supp1}) and (\ref{furth}), we have that
\begin{align}
v \pm \frac{y}{4x} \gg  \Big|1-\frac{y}{4x}\Big|.
\end{align}
By (\ref{Hhat}), we get the bound
\begin{align}
\Big(\frac{1}{x|1-\frac{y}{4x}|} \Big)^B \Big(\frac{x}{K} \Big)^{2B+1}.
\end{align}
\endproof

For very large $x$, we have
\begin{lemma} \label{bigx}
For $x>K^{2-\ep}, 0<y<K^{2+\ep}$ and $h$ as in Lemma \ref{noik}, we have that
\begin{align}
\label{avg1} \sum_{k \equiv 0 \bmod 2} i^k h\Big(\frac{k}{K}\Big) J_{k-1}( x) J_{2k-1}(y) \ll K^{-5/6}.
\end{align}
\end{lemma}
\proof
The idea here is that $J_{k-1}(x)$ is well understood for such large $x$, so that the lemma reduces to understanding the average of one $J$-Bessel function,  $J_{2k-1}(y)$. Define the following functions of $u$ for $u\in (1,2)$:
\begin{align}
z(u)  =  x\left(1-\frac{(uK-1)^2}{x^2}\right)^{1/2}+  (uK-1)\tan^{-1}\left( \left(\frac{uK-1}{x} \right)\left(1-\frac{(uK-1)^2}{x^2}\right)^{-1/2}\right)= x+O(K^{\ep})
\end{align}
and
\begin{align}
\label{smallh2} h_2(u) =  \frac{\sin z(u) - \cos z(u)}{\sqrt{ z(u)-\frac{(k-1)\pi}{2}}}.
\end{align}
Note that for any integer $B\ge 1$, we have that 
\begin{align}
\label{bbb} h_2^{(B)}(u) \ll_B K^{\ep B}.
\end{align}

By Lemma \ref{langer} we have the estimate
\begin{align}
\label{lang}  J_{k-1}(x) &= \sqrt{\frac{2}{\pi(z(\frac{k}{K})-\frac{(k-1)\pi}{2})}}\cos\Big( z\Big(\frac{k}{K}\Big)-\frac{(k-1)\pi}{2}-\frac{\pi}{4}\Big)+ O(k^{-4/3})\\
 \nonumber &= \frac{i^k}{\sqrt{\pi}} h_2\Big( \frac{k}{K} \Big) + O(k^{-4/3}) .
\end{align}
Thus the left hand side of (\ref{avg1}) is bounded by
\begin{align}
\label{aaa} \Big|\sum_{k \equiv 0 \bmod 2} h\Big( \frac{k}{K} \Big)h_2\Big( \frac{k}{K} \Big) J_{2k-1}(y)\Big| + \frac{1}{K^{4/3}}\sum_{k \equiv 0 \bmod 2} h\Big(\frac{k}{K}\Big) |J_{k-1}(y)|.
\end{align}
Using (\ref{smallh2}-\ref{bbb}) and Lemma \ref{noik}, the first term above is $\ll K^{1-\epsilon}$. By (\ref{jbound}), the second term is $\ll K^{-5/6}$.
\endproof

\subsection{Possion summation}

\begin{lemma} \label{poiss} Let $N,c\ge 1$ and $\Psi(x)$ a smooth function, compactly supported on $(1,2)$ and satisfying
\begin{align}
\label{repint} \Big(\frac{N}{c}\Big)^{-B} \Psi^{(B)}(x) \ll_B N^{D-\ep B}.
\end{align}
for some $D$ and any integer $B\ge 0$. Let $S(x)$ be a periodic function satisfying $S(x)=S(x+c)$. We have that
\begin{align}
\label{poisson} \sum_{n\ge 1} S(n)\Psi\Big(\frac{n}{N}\Big) =  \frac{\hat{\Psi}(0) N }{c} \sum_{a \bmod c}  S(a) + O_A(N^{-A})
\end{align}
for any integer $A>0$.
\end{lemma}
\proof
By Poisson summation we have that
\begin{align}
\sum_{n\ge 1} S(n)\Psi\Big(\frac{n}{N}\Big) = \frac{1}{c}\sum_{a \bmod c} S(a) \sum_{m\in \mathbb{Z}}  e\left( \frac{ma}{c}\right) \intt \Psi\Big(\frac{x}{N}\Big)  e\left( \frac{-mx}{c}\right) \ dx.
\end{align}
The term with $m=0$ constitutes the main term of (\ref{poisson}). Integrating by parts repeatedly and using (\ref{repint}) shows that the contribution of the terms with $m\neq 0$ falls into the error term.
\endproof

\section{Main term}

By (\ref{watson}) and the approximate functional equation, we have that
\begin{align}
 \| F \|_4^4 = \frac{2\pi^3}{(2k-1)L(1, \text{sym}^2 f)^2}  \sum_{m,n,r\ge 1} \frac{A_f(n,r) }{(mnr^2)^{1/2}} V_{k,1}(m)V_{k,2}( n r^2 ) \sum_{g \in B_{2k}} \frac{a_g(m)a_g(n)}{L(1, \text{sym}^2 g)}.
 \end{align}
Now by the Petersson trace formula, we get that
 \begin{multline}
\| F \|_4^4 = \frac{\pi}{L(1, \text{sym}^2 f)^2} \Bigg( \sum_{n,r\ge 1} \frac{ A_f(n,r) V_{k,1}(n)V_{k,2}( n r^2)}{nr} \\
 + 2\pi \sum_{m,n,r\ge 1}  \frac{A_f(n,r) V_{k,1}(m)V_{k,2}( n r^2)}{(mnr^2)^{1/2}}  \sum_{c_1\ge 1} \frac{S(n,m;c_1)}{c_1}J_{2k-1}\left(\frac{4\pi\sqrt{nm}}{c_1}\right)\Bigg).
\end{multline}
We have that
\begin{multline}
\label{mainterm} \sum_{n,r\ge 1} \frac{ A_f(n,r) V_{k,1}(n)V_{k,2}( n r^2)}{nr} \\
= \frac{1}{2\pi i}\int_{(2)} \frac{1}{2\pi i} \int_{(2)} \frac{\Lambda_{k,1}(1/2+s_1)}{\Lambda_{k,1}(1/2)}\frac{\Lambda_{k,2}(1/2+s_2)}{\Lambda_{k,2}(1/2)} \frac{L(1+s_1+s_2,\sym f)L(1+2s_2, \sym f)}{\zeta(2+2s_1+2s_2)} \ \frac{ds_1}{s_1}\frac{ds_2}{s_2},
\end{multline}
on using (\ref{vbound}) and  the identity
\begin{align}
\sum_{n,r\ge 1} \frac{A_f(n,r)}{n^{s}r^w} = \frac{L(s,\sym f)L(w,\sym f)}{\zeta(s+w)}
\end{align}
of Bump (cf. \cite[Proposition 6.6.3]{gol}), valid for $\Re(s), \Re(w)>1$. We shift the lines of integration in (\ref{mainterm}) to $\Re(s_1)=\Re(s_2)=-1/4$, pick up residues at $\Re(s_1)=\Re(s_2)=0$ and bound the new integral using (\ref{lambdabound}). The result is that
\begin{align}
 \sum_{n,r\ge 1} \frac{ A_f(n,r) V_{k,1}(n)V_{k,2}( n r^2)}{nr} = \frac{6}{\pi^2} L(1,\sym f)^2 + O(k^{-1/2}).
\end{align}
Thus
\begin{multline}
 \label{maindone} \| F \|_4^4= \frac{6}{\pi} \\
  +\frac{2\pi^2}{L(1, \text{sym}^2 f)^2}  \sum_{m,n,r\ge 1}  \frac{A_f(n,r) V_{k,1}(m)V_{k,2}( n r^2)}{(mnr^2)^{1/2}}  \sum_{c_1\ge 1} \frac{S(n,m;c_1)}{c_1}J_{2k-1}\left(\frac{4\pi\sqrt{nm}}{c_1}\right) + O(k^{-1/2}).
\end{multline}
When (\ref{maindone}) is inserted into (\ref{ffrom}), we see that the constant  $\frac{6}{\pi}$ gives the desired main term of Theorem \ref{thm}, on using Lemma \ref{remove} again. It remains to prove that the contribution of the second line of (\ref{maindone}) falls into the error term.

\section{Error term}

It remains to prove that for some $\delta>0$ and any $1\le d_1,d_2,d_3\le k^\ep$, we have that
\begin{multline}
\label{offdiag} \frac{1}{K^2} \sum_{k \equiv 0 \bmod 2} w\left(\frac{k}{K}\right) \sum_{f \in B_k} \frac{1}{L(1, \text{sym}^2 f)}\\
 \sum_{m,n,r,c\ge 1}  \frac{A_f(n,r) a_f(d_1^2d_2^2)  V_{k,1}(m)V_{k,2}( n r^2) }{(mnr^2)^{1/2}} \frac{S(n,m;c_1)}{c_1}J_{2k-1}\left(\frac{4\pi\sqrt{nm}}{c_1}\right) \ll K^{-\delta}.
\end{multline}
Expressing the $GL(3)$ coefficents in terms of $GL(2)$ coefficients, we get
\begin{align}
& \sum_{n, r\ge 1} \frac{A_f(n,r)V_{k,2}(n r^2)S(n,m;c_1)J_{2k-1}(4\pi \sqrt{nm}/c_1)V_{k,2}(n r^2)}{n^{1/2}r}\\
 = & \sum_{n,r,\alpha \ge 1} \frac{A_f(n,1)A_f(r,1)S(n\alpha,m;c_1)J_{2k-1}(4\pi \sqrt{n\alpha m}/c_1)V_{k,2}(n r^2 \alpha^3 )}{n^{1/2}r\alpha^{3/2}}\\
  = & \sum_{n,r,\alpha, \beta, \gamma \ge 1} \frac{a_f(n^2) a_f(r^2) S(n\alpha\beta^2 ,m;c_1)J_{2k-1}(4\pi \sqrt{n\alpha \beta^2 m}/c_1)V_{k,2}(n r^2 \alpha^3 \beta^2 \gamma^4 )}{n^{1/2}r\alpha^{3/2}\beta\gamma^2}.
\end{align}
Finally we can combine $a_f(r^2)$ and $a_f(d_1^2d_2^2)$ using the Hecke relations (\ref{heckerel}) to see that to establish (\ref{offdiag}), it is suffices to prove
\begin{proposition}\label{prop} For some $\delta>0$ and any integers $\alpha,\beta,\gamma,r_1,r_2\ge 1$ such that $k^{-\ep} < r_1/r_2 < k^{\ep}$, we have that
\begin{multline}
\label{propline} \frac{1}{K^2} \sum_{k \equiv 0 \bmod 2} w\left(\frac{k}{K}\right) \sum_{f \in B_k} \frac{1}{L(1, \text{sym}^2 f)} \sum_{n,m,c_1\ge 1} \frac{a_f(n^2)a_f(r_1^2)}{(mn)^{1/2}} \frac{S(n\alpha \beta^2,m;c_1)}{c_1}\\
\times J_{2k-1}\left(\frac{4\pi\sqrt{nm\alpha\beta^2}}{c_1}\right)V_{k,1}(m)V_{k,2}(nr_2^2\alpha^3 \beta^2 \gamma^2) \ll K^{-\delta}.
\end{multline}
\end{proposition}
\noindent The rest of the section is devoted to proving this. By the Petersson trace formula, we have that 
\begin{align}
(\ref{propline}) = \mathcal{E}_1 + \mathcal{E}_2,
\end{align}
where
\begin{multline}
\label{2pet}  \mathcal{E}_1= \frac{1}{2\pi^2K} \sum_{k \equiv 0 \bmod 2} w\left(\frac{k}{K}\right) \frac{k-1}{K} \sum_{m,c_1\ge 1}  \frac{S(r_1\alpha \beta^2,m;c_1)}{(mr_1)^{1/2}c_1}J_{2k-1}\left(\frac{4\pi\sqrt{r_1 m\alpha \beta^2}}{c_1}\right)\\
\times V_{k,1}(m)V_{k,2}(r_1r_2^2\alpha^3 \beta^2 \gamma^2).
\end{multline}
and
\begin{multline}
\mathcal{E}_2 =  \frac{1}{\pi K} \sum_{k \equiv 0 \bmod 2} i^k w\left(\frac{k}{K}\right) \frac{k-1}{K} \sum_{n,m,c_1,c_2\ge 1} \frac{1}{(nm)^{1/2}} \frac{S(n^2,r_1^2;c_2)}{c_2} \frac{S(n\alpha\beta^2,m;c_1)}{c_1} \\
\times J_{k-1}\left(\frac{4\pi n r_1}{c_2}\right)J_{2k-1}\left(\frac{4\pi\sqrt{n m\alpha \beta^2}}{c_1}\right) V_{k,1}(m)V_{k,2}(n r_2^2\alpha^3 \beta^2 \gamma^2).
\end{multline}
We split $\mathcal{E}_2$ further into three pieces,
\begin{align}
\mathcal{E}_2 = E_1 + E_2 + E_3,
\end{align}
where $E_1$ consists of those terms of $\mathcal{E}_2$ with
\begin{align}
\frac{nr_1}{c_2}> K^{2-\epsilon},
\end{align}
$E_2$ consists of those terms with
\begin{align}
\label{e21} K^{4/3-\epsilon}< \frac{nr_1}{c_2} < K^{2-\ep} \ \text{  and  } \ \frac{nr_1}{c_2 K^{\epsilon}}<\frac{\sqrt{n m\alpha \beta^2}}{c_1}<\frac{nr_1K^{\epsilon}}{c_2 },
\end{align}
and $E_3$ is the rest of $\mathcal{E}_2$.

We first deal with $\mathcal{E}_1$. 
We have
\begin{lemma} 
\begin{align}
\mathcal{E}_1\ll K^{-1/4+\epsilon}.
\end{align}
\end{lemma}
\proof
By (\ref{v-h}), (\ref{klobound}) and Lemma \ref{noik} we have that
\begin{align}
\mathcal{E}_1\ll \frac{1}{K^{1-\epsilon}} \sum_{\substack{m<K^{1+\epsilon}\\ c_1 < K^{1/2+\epsilon}}} \frac{1}{m^{1/2}} \frac{1}{c_1^{1/2}}  \ll K^{-1/4+\epsilon}.
\end{align}
\endproof

As for $\mathcal{E}_2$ we first show that
\begin{lemma} \label{e1lemma}
\begin{align}
 E_1 \ll K^{-1/13}.
\end{align}
\end{lemma}
\proof
Since 
\begin{align}
nr_1/c_2 > K^{2-\epsilon},
\end{align}
we can assume by (\ref{vbound}) and (\ref{jbound}) that $nr_1 <  K^{2+\epsilon}$ and $c_2<K^{\epsilon}$. By (\ref{klobound}), (\ref{v-h}) and Lemma \ref{bigx}, we have that $E_1$ is less than
\begin{align}
\frac{1}{K^{11/6-\ep}}  \sum_{\substack{n< K^{2+\epsilon}\\  m<K^{1+\ep}\\ c_1<K^{1/2+\epsilon} } }  \frac{1}{( nmc_1)^{1/2}} \ll K^{-11/13}
\end{align}
for $\ep$ small enough.
\endproof
For the second piece we have that
\begin{lemma} \label{e2lemma}
\begin{align}
E_2 \ll K^{-1/15}.
\end{align}
\end{lemma}
\proof
By (\ref{vbound}), (\ref{jbound}) and (\ref{e21}), we may assume that the terms in $E_2$ satisfy $n< K^{2+\ep}/r_1^2$, $m<K^{1+\ep}$, $c_2 > \sqrt{n}r_1c_1/ K^{1/2-\ep}\sqrt{\alpha}\beta$, and $c_1< K^{1/6+\ep}$.

For each choice of $n,c_1,c_2$, let $m_0$ be the rational number such that
\begin{align}
\frac{ 4\pi \sqrt{n m_0\alpha \beta^2}}{c_1}=\frac{16\pi n r_1}{c_2}.
\end{align}
Consider first the contribution to $E_2$ of the terms with 
\begin{align}
\label{condi} |m-m_0| < m_0\cdot \min\Big\{\frac{1}{2}, \frac{ n r_1}{c_2 K^{199/100}} \Big\} = \min\Big\{ \frac{m_0}{2},  \frac{ n^2 r_1^3 c_1^2}{c_2^3 \alpha\beta^2 K^{199/100}} \Big\},
\end{align}
so that
\begin{align}
\frac{1}{\sqrt{m}}\ll \frac{1}{\sqrt{m_0}} = \frac{c_2\sqrt{\alpha} \beta}{4c_1 r_1 \sqrt{n}}.
\end{align}
By (\ref{klobound}) and the bound $Kc_2/nr_1$ of (\ref{doublebound1}), we see that the contribution of the terms satisfying (\ref{condi}) is less than
\begin{align}
\frac{1}{K} \sum_{\substack{ n< K^{2+\ep}/r_1^2  \\ c_1< K^{1/6+\ep}\\  c_2 > \sqrt{n}r_1c_1/ K^{1/2-\ep}\sqrt{\alpha} \beta }}
 \frac{ n^2 r_1^3 c_1^2}{c_2^3 \alpha\beta^2 K^{199/100}}  \frac{c_2\sqrt{\alpha}\beta}{c_1 r_1 \sqrt{n}} \frac{1}{\sqrt{n}}  \frac{1}{\sqrt{c_1}}\frac{1}{\sqrt{c_2}} \frac{Kc_2}{nr_1} \ll K^{-11/150+\ep}.
\end{align}

Now consider those terms not satisfying (\ref{condi}). If
\begin{align}
|m-m_0|\ge \frac{m_0}{2},
\end{align}
then
\begin{align}
\label{inci} \Big|\frac{4\pi \sqrt{n m\alpha \beta^2}}{c_1}-\frac{16\pi n r_1}{c_2}\Big| \gg \frac{n r_1}{c_2}. 
\end{align}
By (\ref{doublebound2}), we see that the contribution of such terms is bounded by any negative power of $K$, since $nr_1/c_2<K^{2-\ep}$. Now consider the case
\begin{align}
\frac{ n r_1}{c_2 K^{199/100}} < \frac{1}{2}
\end{align}
and
\begin{align}
|m-m_0|\ge\frac{m_0 n r_1}{c_2 K^{199/100}},
\end{align}
for which we have that
\begin{align}
\Big|\frac{4\pi \sqrt{n m\alpha \beta^2}}{c_1}-\frac{16\pi n r_1}{c_2}\Big|  \gg  \frac{n r_1}{c_2} \frac{nr_1}{c_2K^{199/100}}.
\end{align}
By (\ref{doublebound2}), we see that the contribution of these terms is bounded by
 \begin{align}
  \frac{n r_1}{Kc_2} K^{-B/100}
 \end{align}
 for any integer $B\ge 0$. 
\endproof

It remains to show that $E_3$ is less than a negative power of $K$. 
\begin{lemma}
\begin{align}
E_3 \ll K^{-1/4+\ep}.
\end{align}
\end{lemma}
\proof
We set
\begin{align}
h(u) = w(u)  V_{uK,1}(m)V_{uK,2}(n r_2^2\alpha^3 \beta^2 \gamma^2),
\end{align}
and apply Lemma \ref{doubleavg} to the sum in $E_3$:
\begin{align}
\label{applyingdoubleavg} &\sum \frac{S(n^2,r_1^2;c_2) S(n\alpha\beta^2,m;c_1)}{(nm)^{1/2} c_1 c_2 K} \sum_{k \equiv 0 \bmod 2}  i^k h\left(\frac{k}{K}\right) J_{k-1}\left(\frac{4\pi n r_1}{c_2}\right)J_{2k-1}\left(\frac{4\pi\sqrt{n m\alpha \beta^2}}{c_1}\right)\\
\nonumber =  &\sum  \frac{S(n^2,r_1^2;c_2) S(n\alpha\beta^2,m;c_1)}{(mc_2 r_1)^{1/2} n c_1  K}  \sum_{\pm} e\left( \pm \left(  \frac{mc_2 \alpha\beta^2}{4c_1^2 r_1} +  \frac{2 n r_1}{c_2} \right) \right)\\ 
 \nonumber &\times  \left(1\pm i \frac{\sqrt{m\alpha \beta^2}c_2}{4 c_1 r_1 \sqrt{n} }  \Big/\sqrt{1-\frac{m\alpha\beta^2 c_2^2}{16 n c_1^2 r_1^2 }}\right)
 H^{\pm} \left( \frac{2\pi \sqrt{nm\alpha \beta^2}}{c_1K} \sqrt{1-\frac{m\alpha\beta^2 c_2^2}{16 n c_1^2 r_1^2 }} \right) + O_A(K^{-A}).
\end{align}
We seek cancellation only in the $n$-sum. To study it in dyadic intervals, let $\psi(x)$ be a smooth, compactly supported function on $(1,2)$ with bounded derivatives and let
\begin{align}
\Psi(x) = \frac{ \psi(x)}{x}  \left(1\pm i \frac{\sqrt{m\alpha \beta^2}c_2}{4 c_1 r_1 \sqrt{xN} }  \Big/\sqrt{1-\frac{m\alpha\beta^2 c_2^2}{16 xN c_1^2 r_1^2 }}\right) H^{\pm} \left( \frac{2\pi \sqrt{x Nm\alpha \beta^2}}{c_1K} \sqrt{1-\frac{m\alpha\beta^2 c_2^2}{16 x N c_1^2 r_1^2 }}  \right)
\end{align}
for $K^{1-\ep}<N<K^{2+\ep}/r_2^2\alpha^{3}\beta^{2}\gamma^{2}$.
Using (\ref{pro1}), we have that
\begin{align}
\Psi^{(B)}(x) \ll_B \left(\frac{\sqrt{Nm\alpha \beta^2}}{c_1K}\right)^{2B+1}
\end{align}
for any integer $B\ge 0$. 
Lest $\Psi(x)=0$, we may assume by the observation made at the end of Lemma \ref{doubleavg} that
\begin{align}
c_2 \ll \frac{Nr_1}{K}.
\end{align}
We may also assume that
\begin{align}
\label{yboun} \frac{\sqrt{N m\alpha \beta^2}}{c_1}<K^{4/3-\ep},
\end{align}
because if not then $\Psi(x)$ vanishes unless 
\begin{align}
\label{donealready} \frac{Nr_1}{c_2} \sim \frac{\sqrt{N m\alpha \beta^2}}{c_1}.
\end{align}
But the terms satisfying (\ref{donealready}) with $ \sqrt{N m\alpha \beta^2}/c_1>K^{4/3-\ep}$ belong to $E_2$ and were treated in Lemma \ref{e2lemma}.
Thus by (\ref{vbound}) we may assume that
\begin{align}
\label{extra} \left(\frac{c_1c_2}{N}\frac{Nm\alpha \beta^2}{c_1^2K^2 }\right)^B \ll_B \left(\frac{Nm\alpha \beta^2 r_1}{c_1 K^3}\right)^B <K^{-B/6}
\end{align}
for any integer $B\ge 0$.
Applying Lemma \ref{poiss} we get that
\begin{multline}
\sum_{n\ge 1} \frac{S(n^2,r_1^2;c_2) S(n\alpha\beta^2,m;c_1)}{N}  e\left( \frac{\pm 2 n r_1}{c_2} \right) \Psi\Big(\frac{n}{N}\Big)\\
= \frac{\hat{\Psi}(0)}{c_1c_2} \sum_{a\bmod c_1c_2} S(a^2,r_1^2;c_2)e\left( \frac{2 a r_1}{c_2} \right) S(a\alpha \beta^2, m ; c_1) + O_A(K^{-A}).
\end{multline}
Bounding the rest of the sum in $E_3$ absolutely, we get that
\begin{align}
 \label{afterpoiss}E_3 \ll \frac{1}{K^{1-\ep}} \sum_{\substack{m<K^{1+\ep}\\ c_1<K^{1/2+\ep}\\c_2< \frac{K^{1+\ep}c_1}{m^{1/2}\alpha\beta^2}}}  \frac{1}{m^{1/2} {c_2}^{3/2}c_1^2} \left| \hat{\psi}(0)\sum_{a\bmod c_1c_2} S(a^2,r_1^2;c_2)e\left( \frac{2 a r_1}{c_2} \right) S(a\alpha \beta^2, m ; c_1) \right|,
 \end{align}
where the range of summation in $m,c_1,c_2$ is gotten from (\ref{vbound}) and the fact that (\ref{Hpm}) is nonzero only for $y<4x$.. By (\ref{pro2}), we have that $\hat{\psi}(0) \ll 1$.

Write $c_1=b_1 c_1'$ and $c_2= b_2 c_2'$, where $(b_1,c_1')=(b_2,c_2')=(c_1',c_2')=1$ and $p|b_1 \Leftrightarrow p|b_2$. By (\ref{kloosmult}), we have that (\ref{afterpoiss}) is bounded by
\begin{align}
\label {threesums} \frac{1}{K^{1-\ep}} \sum_{\substack{b_1,b_2\le K^{10} \\ p|b_1\Leftrightarrow p|b_2}} \sum_{\substack{m<K^{1+\ep}\\ c_1'<K^{1/2+\ep} \\c_2'< \frac{K^{1+\ep}c_1'}{b_2 m^{1/2}\alpha \beta^2}}} \frac{1}{m^{1/2}}|S_1 S_2 S_3|,
\end{align}
where
\begin{align}
&S_1 = \frac{1}{{c_2'}^{3/2}} \sum_{a \bmod c_2'} S(a^2, r_1^2 \overline{{b_2}}^2; c_2')e\left( \frac{2 a r_1 \overline{b_2}}{c_2'} \right),\\
&S_2 =  \frac{1}{{c_1'}^{2}} \sum_{a \bmod c_1'} S(a \alpha \beta^2, m ; c_1'),\\
&S_3=  \frac{1}{ {b_2}^{3/2}b_1^2} \sum_{a\bmod b_1b_2} S(a^2,r_1^2\overline{c_2'}^2;b_2)e\left( \frac{2 a r_1 \overline{c_2'} }{b_2} \right) S(a\alpha \beta^2, m \overline{c_1'}^2 ; b_1).
 \end{align}
Now we provide bounds for these sums. By \cite[Lemma 3.3]{kha}, we have that
\begin{align}
S_1 \le  \delta_{c_2' = \square},
\end{align}
while it is easy to see that
\begin{align}
S_2 \le  \delta_{c_1' | \alpha \beta^2},
\end{align}
and by Weil's bound we have that
\begin{align}
S_3\ll  \frac{K^{\ep}}{ {b_2}b_1^{3/2}} \sum_{1\le a \le b_1b_2}  (a^2,r_1^2,b_2)^{1/2} (a\alpha \beta^2, m , b_1)^{1/2} \ll   \frac{ K^{\ep}(m , b_1)^{1/2}}{ {b_2} b_1^{3/2}} \sum_{1\le a \le b_1b_2}  (a,r_1) \ll   \frac{ K^{\ep} (m , b_1)^{1/2}}{ b_1^{1/2}}.
\end{align}
Thus (\ref{threesums}) is bounded by
\begin{align}
&\frac{1}{K^{1-\ep}}  \sum_{\substack{b_1,b_2\le K^{10} \\ p|b_1\Leftrightarrow p|b_2}}  \frac{1}{b_1^{1/2}} \sum_{m\le K^{1+\ep}} \frac{ (m , b_1)^{1/2}}{ m^{1/2}} \sum_{c_1' | \alpha \beta^2}  \sum_{c_2'< \frac{K^{1+\ep}c_1' }{b_2 m^{1/2}\alpha \beta^2}}  \delta_{c_2' = \square}\\
\nonumber \ll &\frac{1}{K^{1/2-\ep}} \sum_{\substack{b_1,b_2\le K^{10} \\ p|b_1\Leftrightarrow p|b_2}}  \frac{1}{b_1^{1/2} b_2^{1/2}} \sum_{m\le K^{1+\ep}} \frac{ (m , b_1)^{1/2}}{ m^{3/4}}\\
\nonumber \ll &\frac{1}{K^{1/4-\ep}} \sum_{\substack{b_1,b_2\le K^{10} \\ p|b_1\Leftrightarrow p|b_2}}  \frac{1}{b_1^{1/2} b_2^{1/2}}\\
\nonumber \ll &\frac{1}{K^{1/4-\ep}}.
\end{align}
\endproof

\bibliographystyle{amsplain}

\bibliography{averages4}

\providecommand{\bysame}{\leavevmode\hbox to3em{\hrulefill}\thinspace}
\providecommand{\MR}{\relax\ifhmode\unskip\space\fi MR }
\providecommand{\MRhref}[2]{%
  \href{http://www.ams.org/mathscinet-getitem?mr=#1}{#2}
}
\providecommand{\href}[2]{#2}
\begin{thebibliography}{10}

\bibitem{ber}
Michael Berry, \emph{Regular and irregular semiclassical wavefunctions}, J.
  Phys. A \textbf{10} (1977), no.~12, 2083--2091.

\bibitem{blokhayou}
Valentin Blomer, Rizwanur Khan, and Matthew Young, \emph{Mass distribution of
  holomorphic cusp forms}, preprint.

\bibitem{erd}
Arthur Erd{\'e}lyi, Wilhelm Magnus, Fritz Oberhettinger, and Francesco~G.
  Tricomi, \emph{Higher transcendental functions. {V}ol. {II}}, Robert E.
  Krieger Publishing Co. Inc., Melbourne, Fla., 1981, Based on notes left by
  Harry Bateman, Reprint of the 1953 original.

\bibitem{gol}
Dorian Goldfeld, \emph{Automorphic forms and {$L$}-functions for the group
  {${\rm GL}(n,\bold R)$}}, Cambridge Studies in Advanced Mathematics, vol.~99,
  Cambridge University Press, Cambridge, 2006, With an appendix by Kevin A.
  Broughan.

\bibitem{hejrac}
Dennis Hejhal and Barry Rackner, \emph{On the topography of {M}aass waveforms
  for {${\rm PSL}(2,{\bf Z})$}}, Experiment. Math. \textbf{1} (1992), no.~4,
  275--305.

\bibitem{holsou}
Roman Holowinsky and Kannan Soundararajan, \emph{Mass equidistribution for
  {H}ecke eigenforms}, Ann. of Math. (2) \textbf{172} (2010), no.~2,
  1517--1528.

\bibitem{iwa}
Henryk Iwaniec, \emph{Topics in classical automorphic forms}, Graduate Studies
  in Mathematics, vol.~17, American Mathematical Society, Providence, RI, 1997.

\bibitem{iwakow}
Henryk Iwaniec and Emmanuel Kowalski, \emph{Analytic number theory}, American
  Mathematical Society Colloquium Publications, vol.~53, American Mathematical
  Society, Providence, RI, 2004.

\bibitem{kha}
Rizwanur Khan, \emph{Non-vanishing of the symmetric square {$L$}-function at
  the central point}, Proc. Lond. Math. Soc. (3) \textbf{100} (2010), no.~3,
  736--762.

\bibitem{lauwu}
Yuk-Kam Lau and Jie Wu, \emph{A density theorem on automorphic {$L$}-functions
  and some applications}, Trans. Amer. Math. Soc. \textbf{358} (2006), no.~1,
  441--472 (electronic).

\bibitem{liuyou}
Sheng-Chi Liu and Matthew Young, \emph{Growth and nonvanishing of restricted
  {S}iegel modular forms arising as {S}aito-{K}urokawa lifts}, preprint.

\bibitem{ran}
Robert Rankin, \emph{The vanishing of {P}oincar\'e series}, Proc. Edin. Math.
  Soc. \textbf{23} (1980), no.~2, 151--161.

\bibitem{wat}
Thomas Watson, \emph{Rankin triple products and quantum chaos}, to appear in
  Ann. of Math. (2).

\end{thebibliography}

\end{document}